\documentclass[11pt]{article}

\usepackage[utf8]{inputenc}
\usepackage{amsmath}
\usepackage{amsfonts,amssymb,latexsym,multicol,color}
\usepackage[francais]{babel}
\usepackage[T1]{fontenc}

\newcounter{fig}
\newtheorem{theo}{Th\'eor\`eme}

\newtheorem{prop}{Proposition}

\newcommand{\virg}{\raisebox{.7mm}{,}}

\newcommand{\expli}[1]{\quad\text{\footnotesize (#1)}}

\makeatletter
\ifcase\@ptsize

\or

\or

\fi
\makeatother

\newcommand{\fhi}{\varphi}
\newcommand{\ioe}{\leqslant}
\newcommand{\soe}{\geqslant}

\newcommand{\vers}{\rightarrow}

\newcommand{\Nat}{{\mathbb N}}

\newcommand{\fin}{\hfill$\Box$}
\newcommand{\dem}{\noindent {\bf D\'emonstration\ }}
\newcommand{\fine}{\tag*{\mbox{$\Box$}}}

\providecommand{\bysame}{\leavevmode ---\ }
\providecommand{\og}{``}
\providecommand{\fg}{''}
\providecommand{\smfandname}{et}
\providecommand{\smfedsname}{\'eds.}
\providecommand{\smfedname}{\'ed.}
\providecommand{\smfmastersthesisname}{M\'emoire}
\providecommand{\smfphdthesisname}{Th\`ese}

\title{Sur la variation quadratique totale de la suite des parties fractionnaires des quotients d'un nombre r\'eel positif par les nombres entiers naturels cons\'ecutifs}

\author{Michel Balazard}
\date{}
\begin{document}
\maketitle

\begin{center}
  {\sc Abstract}
\end{center}
\begin{quote}
{\footnotesize We give an asymptotic formula for the quadratic total variation of the sequence of fractional parts of the quotients of a positive real number by the consecutive natural numbers :
$$
\sum_{n\soe 1}\big (\{x/(n+1)\}-\{x/n\}\big )^2=\frac{\zeta(3/2)}{\pi}x^{1/2}+O(x^{3/7}).
$$}
\end{quote}

\begin{center}
  {\sc Keywords}
\end{center}
\begin{quote}
{\footnotesize Fractional part, quadratic total variation\\MSC classification : 11N37}
\end{quote}



\section{\'Enonc\'e du r\'esultat}

Posons
$$
Q(x)=\sum_{n\soe 1}\big (\{x/(n+1)\}-\{x/n\}\big )^2\quad (x>0),
$$
o\`u
$\{t\}=t-\lfloor t\rfloor$ d\'esigne la partie fractionnaire du nombre r\'eel $t$, et $\lfloor t\rfloor$ sa partie enti\`ere. Le pr\'esent document est annexe \`a l'article \cite{bw2015} et a pour objet la d\'emonstration du r\'esultat suivant, mentionn\'e dans \cite{bw2015}.
\begin{theo}
On a
$$
Q(x)=\frac{\zeta(3/2)}{\pi}x^{1/2}+O(x^{3/7}) \quad (x>0).
$$
\end{theo}

\smallskip

La d\'emonstration est similaire \`a celle du th\'eor\`eme de \cite{bw2015}. Nous reprenons l'ensemble des notations de cet article.

\section{D\'emonstration du th\'eor\`eme}

Posons donc, pour $d\in \Nat$,
$$
Q_d(x)=\sum_{n\soe 1}[\lfloor x/n\rfloor -\lfloor x/(n+1)\rfloor =d] \cdot \big(\{x/n\}-\{x/(n+1)\}\big)^2\, ,
$$
o\`u nous avons utilis\'e la notation d'Iverson : $[P]=1$ si la proposition $P$ est vraie, $0$ sinon.

Nous allons \'evaluer la quantit\'e $Q_d(x)$ pour tout $d$.

\subsection{Contribution des grandes valeurs de $d$}

On a d'abord, comme dans le cas de la fonction $W$ de \cite{bw2015} :
\begin{equation}\label{t48}
\sum_{d>D}Q_d(x) <\sqrt{x/(D-1)}
\end{equation}
pour $D>1$ et $x>0$.

\subsection{Estimation de $Q_0(x)$}

On a
\begin{align*}
Q_0(x) &=\sum_{0\ioe k\ioe x}\; \sum_{x/(k+1)< n \ioe x/k-1}x^2/n^2(n+1)^2\\
&=\Delta_1(x)-\Delta_2(x),
\end{align*}
o\`u
\begin{align*}
\Delta_1(x)&=x^2 \sum_{0\ioe k\ioe K}\; \sum_{x/(k+1)< n \ioe x/k}1/n^2(n+1)^2\\
\Delta_2(x)&= x^2\sum_{0\ioe k\ioe K}\, 1/\lfloor x/k\rfloor^2(\lfloor x/k\rfloor+1)^2 \, ,
\end{align*}
et
$$
K=K(x)=\lfloor \sqrt{x+1/4}-1/2\rfloor .
$$

Si $x\soe 2$, on a $K\soe 1$ et
\begin{align*}
\Delta_1(x)&=x^2 \sum_{n>x/(K+1)}\big (n^{-4}+O(n^{-5})\big)\\
&=x^2\Big(\frac{K^3+O(K^2)}{3x^3}+O(K^4/x^4)\Big)\\
&=\frac{\sqrt{x}}{3}+O(1)\, ,
\end{align*}
et
\begin{align*}
\Delta_2(x)&= x^2\sum_{0\ioe k\ioe K}\big(k^4/x^4+O(k^5/x^5)\big) \\
&=\frac{K^5+O(K^4)}{5x^2}+O(K^6/x^3)\\
&=\frac{\sqrt{x}}{5}+O(1).
\end{align*}
Par cons\'equent,
\begin{equation}\label{t49}
Q_0(x)=\frac{2\sqrt{x}}{15}+O(1)
\end{equation}

\subsection{Calcul de $Q_d(x)$ pour $d$ positif}

Avec des notations correspondant \`a celles de \cite{bw2015}, nous \'evaluons maintenant
\begin{align*}
Q_{d,1}(x)&= \sum_{K_{d-1}<k\ioe K_d-1}\, \sum_{x/(k-d+1)-1< n \ioe x/k}\Big( d-\frac{x}{n(n+1)}\Big)^2\\
&=\sum_{K_{d-1}<k\ioe K_d-1}\,\bigg(d^2\big(\lfloor x/k\rfloor -\lfloor x/(k-d+1)\rfloor +1\big)-2dx\Big(\frac{1}{\lfloor x/(k-d+1)\rfloor}-\frac{1}{\lfloor x/k\rfloor +1}\Big)\\
&\qquad \qquad +x^2\big(F(x/(k-d+1))-F(x/k+1)\big)\bigg)\, ,
\end{align*}
o\`u
$$
F(t)=\sum_{n>t-1}\frac{1}{n^2(n+1)^2}\cdotp
$$

On se ram\`ene ensuite \`a des sommes $(d-1)$-t\'elescopiques en utilisant la fonction
$$
\fhi(t)=\frac{1}{\lfloor t \rfloor}-\frac{1}{\lfloor t \rfloor +1} \virg
$$
et en notant que
$$
F(t+1)=F(t)-\fhi(t)^2.
$$

Ainsi,
\begin{multline*}
Q_{d,1}(x)=d^2\sum_{K_{d}-d<k\ioe K_d-1}\lfloor x/k\rfloor -d^2\sum_{K_{d-1}-d+1<k\ioe K_{d-1}}\lfloor x/k\rfloor +d^2(K_{d}-K_{d-1}-1) \\
+2dx\sum_{K_{d}-d<k\ioe K_d-1}\frac{1}{\lfloor x/k\rfloor} -2dx \sum_{K_{d-1}-d+1<k\ioe K_{d-1}}\frac{1}{\lfloor x/k\rfloor} -2dx\sum_{K_{d-1}<k\ioe K_d-1}\fhi(x/k) \\
-x^2\sum_{K_{d}-d<k\ioe K_d-1}F(x/k) + x^2\sum_{K_{d-1}-d+1<k\ioe K_{d-1}}F(x/k) +x^2\sum_{K_{d-1}<k\ioe K_d-1}\fhi(x/k)^2.
\end{multline*}

\smallskip

De m\^eme,
\begin{multline*}
Q_{d,2}(x)=-d^2\sum_{K_{d+1}-d-1<k\ioe K_{d+1}}\lfloor x/k\rfloor +d^2\sum_{K_{d}-d<k\ioe K_{d}+1}\lfloor x/k\rfloor -d^2(K_{d+1}-K_{d}-1) \\
-2dx\sum_{K_{d+1}-d-1<k\ioe K_{d+1}}\frac{1}{\lfloor x/k\rfloor} +2dx \sum_{K_{d}-d<k\ioe K_{d}+1}\frac{1}{\lfloor x/k\rfloor} +2dx\sum_{K_{d}+1<k\ioe K_{d+1}}\fhi(x/k) \\
+x^2\sum_{K_{d+1}-d-1<k\ioe K_{d+1}}F(x/k) -x^2 \sum_{K_{d}-d<k\ioe K_{d}+1}F(x/k) -x^2\sum_{K_{d}+1<k\ioe K_{d+1}}\fhi(x/k)^2.
\end{multline*}

\smallskip

Enfin,
\begin{multline*}
Q_{d,3}(x)= \sum_{\frac{x}{K_d+1}<k\ioe \frac{x}{K_d}}\, \Big( d-\frac{x}{n(n+1)}\Big)^2\\
=d^2\big(\lfloor x/K_d \rfloor - \lfloor x/(K_d+1) \rfloor\big)-2dx\Big(\frac{1}{\lfloor x/(K_d+1) \rfloor +1}-\frac{1}{\lfloor x/K_d \rfloor +1}\Big) \\
+x^2\Big( F\big(x/(K_d+1) +1\big)-F(x/K_d+1)\Big).
\end{multline*}

L'addition des formules obtenues pour $Q_{d,1}$, $Q_{d,2}$ et $Q_{d,3}$ donne
\begin{multline}\label{t56}
Q_{d}(x)=\\
\Big(2\sum_{K_{d}-d<k\ioe K_{d}}-\sum_{K_{d+1}-d-1<k\ioe K_{d+1}}-\sum_{K_{d-1}-d+1<k\ioe K_{d-1}}\Big)\big(d^2\lfloor x/k\rfloor +2dx/\lfloor x/k\rfloor-x^2F(x/k)\big) \\
-\Big(\sum_{K_{d} < k\ioe K_{d+1}}-\sum_{K_{d-1}<k\ioe K_{d}}\Big)\big(d-x\fhi(x/k)\big)^2.
\end{multline}

\subsection{Estimation de $Q_d(x)$ pour $d$ positif}

\subsubsection{Estimations compl\'ementaires concernant $K_d$}

Les propositions suivantes n'\'etaient pas utiles pour d\'emontrer le r\'esultat de \cite{bw2015}, mais nous en aurons l'usage pour l'estimation de $Q_d(x)$.
\begin{prop}\label{t62}
Pour $0\ioe d \ioe x$, on a
$$
K_d(x)=\sqrt{dx}+\frac d2 + O(d^{3/2}x^{-1/2})+O(1).
$$
\end{prop}
\dem

On a
\begin{align*}
K_d(x) &= \lfloor (d+\sqrt{d^2+4dx}\,)/2\rfloor\\
&=\sqrt{dx}\,(1+d/4x)^{1/2} +d/2+O(1)\\
&=\sqrt{dx}+\frac d2 + O(d^{3/2}x^{-1/2})+O(1).\fine
\end{align*}

\begin{prop}\label{t69}
Pour $0\ioe d \ioe x$, on a
$$
\sum_{K_{d}-d<k\ioe K_{d}}k=d\sqrt{dx} +O(d^{5/2}x^{-1/2})+O(d).
$$
\end{prop}
\dem

On a
\begin{align*}
\sum_{K_{d}-d<k\ioe K_{d}}k&=\frac 12 \big(K_d^2+K_d-(K_{d}-d)^2-K_{d}+d\big)\\
&=dK_d-\frac{d^2}2+\frac{d}2\\
&=d\sqrt{dx} +O(d^{5/2}x^{-1/2})+O(d)\, ,
\end{align*}
d'apr\`es la proposition \ref{t62}.\fin

\begin{prop}\label{t70}
Pour $0\ioe d \ioe x$, on a
$$
\sum_{K_{d}-d<k\ioe K_{d}}k^3=x\big(d^2\sqrt{dx}+O(d^{7/2}x^{-1/2})+O(d^2)\big).
$$
\end{prop}
\dem

On a
\begin{align*}
\sum_{K_{d}-d<k\ioe K_{d}}k^3&=\frac 14 \big(K_d^4-(K_{d}-d)^4\big)+\frac 12 \big(K_d^3-(K_{d}-d)^3\big)+\frac 14 \big(K_d^2-(K_{d}-d)^2\big)\\
&=dK_d^3-\frac 32 d(d-1)K_d^2 +O(d^3K_d)\expli{puisque $d\ioe K_d$}\\
&=d\big(d^{3/2}x^{3/2}+\frac 32 d^2x +O(d^{5/2}x^{1/2})+O(dx)\big) -\frac 32 d(d-1)\big(dx+O(d^{3/2}x^{1/2})\big)+O(d^{7/2}x^{1/2})\\
&\expli{d'apr\`es la proposition \ref{t62}}\\
&=x\big(d^2\sqrt{dx}+O(d^{7/2}x^{-1/2})+O(d^2)\big).\fine
\end{align*}

\begin{prop}\label{t68}
Pour $0\ioe d \ioe x$, on a
$$
\sum_{K_d<k\ioe K_{d+1}}k^4=\frac{(d+1)^2\sqrt{d+1}-d^2\sqrt{d}}5 x^{5/2}+O((d+1)^2x^2).
$$
\end{prop}
\dem

Soit $P$ le polyn\^ome tel que 
$$
\sum_{k\ioe K}k^4=P(K) \quad (K\in \Nat).
$$

Le terme de plus haut degr\'e de $P(K)$ est $K^5/5$. Par cons\'equent, si $d$ est fix\'e, on a
$$
\sum_{K_d<k\ioe K_{d+1}}k^4 \sim \frac{(d+1)^2\sqrt{d+1}-d^2\sqrt{d}}5 x^{5/2} \quad (x \vers \infty).
$$

Maintenant, soit 
$$
Q(X,Y)=\frac{P(X)-P(Y)}{X-Y}\cdotp
$$
Cette fraction rationnelle est en fait un polyn\^ome de degr\'e $4$. Notons $Q_4(X,Y)$ sa partie homog\`ene de degr\'e $4$, et $R=Q-Q_4$, de sorte que $R$ est de degr\'e $3$.

La contribution de chaque mon\^ome de $Q_4$ \`a la quantit\'e $Q_4(K_{d+1},K_d)$ est de la forme $c(d)x^{2}+O\big((d+1)^{5/2}x^{3/2}\big)$ o\`u la fonction $c(d)$ d\'epend du mon\^ome consid\'er\'e, mais est toujours $O\big((d+1)^2\big)$. D'autre part, $R(K_{d+1},K_d)=O\big((d+1)^{3/2}x^{3/2}\big)$. On en d\'eduit que
$$
P(K_{d+1})-P(K_d)=(K_{d+1}-K_d)\big(C(d)x^{2}+O((d+1)^{5/2}x^{3/2})\big)\, ,
$$
o\`u $C(d)$ est la somme des $c(d)$. En particulier, $C(d)=O\big((d+1)^2\big)$. En utilisant maintenant la proposition 2 de \cite{bw2015}, on obtient
$$
P(K_{d+1})-P(K_d)=C_1(d)x^{5/2}+O((d+1)^{2}x^{2}\big)\, ,
$$
pour une certaine constante $C_1(d)$, qui vaut n\'ecessairement
\begin{equation*}
\frac{(d+1)^2\sqrt{d+1}-d^2\sqrt{d}}5\cdotp\fine
\end{equation*}

\subsubsection{Estimation de la fonction $F$}

Nous utiliserons l'estimation suivante.
\begin{prop}\label{t58}
Pour $t\soe 1$ on a
$$
F(t)=\frac 1{3\lfloor t\rfloor^3} +O(1/t^5).
$$
\end{prop}
\dem

On a l'identit\'e
$$
\frac{3}{n^2(n+1)^2}=\frac 1{n^3}-\frac 1{(n+1)^3} -\frac 1{n^3(n+1)^3}\, ,
$$
donc
\begin{align*}
F(t)&=\sum_{n>t-1}\frac{1}{n^2(n+1)^2}\\
&=\frac 13\sum_{n>t-1}\Big (\frac 1{n^3}-\frac 1{(n+1)^3}\Big) + \sum_{n>t-1}O(n^{-6})\\
&=\frac 1{3\lfloor t\rfloor^3} +O(1/t^5).\fine
\end{align*}

\subsubsection{Estimation pr\'eliminaire concernant les sommants de \eqref{t56}}

La quantit\'e 
$$
d^2\lfloor x/k\rfloor + 2dx/\lfloor x/k\rfloor -x^2F(x/k)
$$
apparaît dans \eqref{t56}. Nous en donnons maintenant une estimation.
\begin{prop}\label{t61}
Pour $0<d\ioe x/2$, on a
$$
d^2\lfloor x/k\rfloor +2dx/\lfloor x/k\rfloor -x^2F(x/k)=2d\sqrt{dx}+kd-\frac{k^3}{3x}+O(d^{5/2}x^{-1/2})\quad (K_d-d<k\ioe K_d).
$$
\end{prop}
\dem

Supposons donc $0<d\ioe x/2$ et $K_d-d<k\ioe K_d$. On a d'abord, d'apr\`es la proposition \ref{t58},
\begin{equation*}\label{t59}
x^2F(x/k)=x^2/3\lfloor x/k\rfloor^3+O(k^5/x^3)\, ,
\end{equation*}
et le terme d'erreur est bien $O(d^{5/2}x^{-1/2})$ puisque $k\ioe K_d \ll \sqrt{dx}$.

Ensuite, la proposition 4 de \cite{bw2015} nous donne
\begin{equation*}\label{t60}
d^2\lfloor x/k\rfloor +dx/\lfloor x/k\rfloor =2d\sqrt{dx} +O(d^{5/2}x^{-1/2}).
\end{equation*}

On a donc
$$
d^2\lfloor x/k\rfloor +2dx/\lfloor x/k\rfloor -x^2F(x/k)=2d\sqrt{dx} +dx/\lfloor x/k\rfloor-x^2/3\lfloor x/k\rfloor^3+O(d^{5/2}x^{-1/2}).
$$

\smallskip

Maintenant
\begin{align*}
dx/\lfloor x/k\rfloor-x^2/3\lfloor x/k\rfloor^3 &=\frac{dx}{x/k-\{x/k\}}-\frac{x^2}{3(x/k-\{x/k\})^3}\\
&=kd\big(1+k\{x/k\}/x+O(k^2/x^2)\big)-\frac{k^3}{3x}\big(1+3k\{x/k\}/x+O(k^2/x^2)\big)\\
&=kd-\frac{k^3}{3x}-(k^4/x^2-k^2d/x)\{x/k\}+O(d^{5/2}x^{-1/2}).
\end{align*}

Enfin
\begin{align*}
k^4/x^2-k^2d/x &=k^2(k+\sqrt{dx})(k-\sqrt{dx})/x^2\\
&\ll d^{5/2}x^{-1/2}
\end{align*}
puisque $k-\sqrt{dx} \ll d$ (cf. \cite{bw2015}, (17)).\fin

\subsubsection{Estimation des trois premi\`eres sommes de \eqref{t56}}

\begin{prop}\label{t63}
Pour $0<d\ioe x/2$, on a
\begin{equation}\label{t71}
\sum_{K_{d}-d<k\ioe K_{d}}\big(d^2\lfloor x/k\rfloor + 2dx/\lfloor x/k\rfloor -x^2F(x/k)\big)=\frac 83 d^2\sqrt{dx} +O(d^{7/2}x^{-1/2})+O(d^2).
\end{equation}
\end{prop}
\dem

La proposition \ref{t61} nous donne :
\begin{equation*}
\sum_{K_{d}-d<k\ioe K_{d}}\big(d^2\lfloor x/k\rfloor + 2dx/\lfloor x/k\rfloor -x^2F(x/k)\big)=2d^2\sqrt{dx} +\sum_{K_{d}-d<k\ioe K_{d}}(kd-k^3/3x)+O(d^{7/2}x^{-1/2})
\end{equation*}

La proposition d\'ecoule de cette expression et des propositions \ref{t69} et \ref{t70}..\fin

\begin{prop}
Pour $0<d+1\ioe x/2$, on a
\begin{multline}\label{t72}
\sum_{K_{d+1}-d-1<k\ioe K_{d+1}}\big(d^2\lfloor x/k\rfloor + 2dx/\lfloor x/k\rfloor -x^2F(x/k)\big)=\\
\frac{8d^2+4d-1}3 \sqrt{(d+1)x} +O(d^{7/2}x^{-1/2})+O(d^2).
\end{multline}
\end{prop}
\dem

On a
\begin{multline}\label{t64}
\sum_{K_{d+1}-d-1<k\ioe K_{d+1}}\big(d^2\lfloor x/k\rfloor + 2dx/\lfloor x/k\rfloor -x^2F(x/k)\big)=\\
\sum_{K_{d+1}-d-1<k\ioe K_{d+1}}\big((d+1)^2\lfloor x/k\rfloor + 2(d+1)x/\lfloor x/k\rfloor -x^2F(x/k)\big)-\\
\sum_{K_{d+1}-d-1<k\ioe K_{d+1}}\big((2d+1)\lfloor x/k\rfloor + 2x/\lfloor x/k\rfloor\big).
\end{multline}

La premi\`ere somme du second membre de \eqref{t64} vaut
$$
\frac 83 (d+1)^2\sqrt{(d+1)x} +O(d^{7/2}x^{-1/2})+O(d^2)
$$
d'apr\`es la proposition \ref{t63}.

En utilisant l'estimation
$$
\frac xk =\sqrt{x/(d+1)}+O(1) \quad (K_{d+1}-d-1<k\ioe K_{d+1})\, ,
$$
on voit que la seconde somme du second membre de \eqref{t64} vaut
$$
\sum_{K_{d+1}-d-1<k\ioe K_{d+1}}\big((2d+1)\lfloor x/k\rfloor + 2x/\lfloor x/k\rfloor\big)=(4d+3)\sqrt{(d+1)x}+O(d^2).
$$

La proposition s'en d\'eduit.\fin

\begin{prop}
Pour $0\ioe d-1\ioe x/2$, on a
\begin{equation}\label{t73}
\sum_{K_{d-1}-d+1<k\ioe K_{d-1}}\big(d^2\lfloor x/k\rfloor + 2dx/\lfloor x/k\rfloor -x^2F(x/k)\big)=\frac{8d^2-4d-1}3 \sqrt{(d-1)x} +O(d^{7/2}x^{-1/2})+O(d^2).
\end{equation}
\end{prop}
\dem

Si $d=1$ le r\'esultat est trivial. Si $d>1$, on a
\begin{multline}\label{t65}
\sum_{K_{d-1}-d+1<k\ioe K_{d-1}}\big(d^2\lfloor x/k\rfloor + 2dx/\lfloor x/k\rfloor -x^2F(x/k)\big)=\\
\sum_{K_{d-1}-d+1<k\ioe K_{d-1}}\big((d-1)^2\lfloor x/k\rfloor + 2(d-1)x/\lfloor x/k\rfloor -x^2F(x/k)\big)\\
+\sum_{K_{d-1}-d+1<k\ioe K_{d-1}}\big((2d-1)\lfloor x/k\rfloor + 2x/\lfloor x/k\rfloor\big).
\end{multline}

La premi\`ere somme du second membre de \eqref{t65} vaut
$$
\frac 83 (d-1)^2\sqrt{(d-1)x} +O(d^{7/2}x^{-1/2})+O(d^2)
$$
d'apr\`es la proposition \ref{t63}.

En utilisant l'estimation
$$
\frac xk =\sqrt{x/(d-1)}+O(1) \quad (K_{d-1}-d+1<k\ioe K_{d-1})\, ,
$$
on voit que la seconde somme du second membre de \eqref{t65} vaut
$$
\sum_{K_{d-1}-d+1<k\ioe K_{d-1}}\big((2d-1)\lfloor x/k\rfloor + 2x/\lfloor x/k\rfloor\big)=(4d-3)\sqrt{(d-1)x}+O(d^2).
$$

La proposition s'en d\'eduit.\fin

\subsubsection{Estimation des deux derni\`eres sommes de \eqref{t56}}

En utilisant la proposition 2 de \cite{bw2015}, on obtient
\begin{equation}\label{t52}
-d^2(K_{d+1}-2K_d+K_{d-1})=-d^2(\sqrt{d+1}+\sqrt{d-1}-2\sqrt{d}\,)\sqrt{x}+O(d^2)
\end{equation}

\smallskip

En utilisant la relation asymptotique
$$
\fhi(t)=\frac 1{t^2} +O(1/t^3) \quad (t\soe 1)
$$
et la proposition 3 de \cite{bw2015}, on obtient
\begin{multline}\label{t53}
2dx\Big (\sum_{K_{d} < k\ioe K_{d+1}}-\sum_{K_{d-1}<k\ioe K_{d}}\Big)\fhi(x/k)=\\
\frac{2d}3 \big((d+1)\sqrt{d+1}+(d-1)\sqrt{d-1}-2d\sqrt{d}\,\big)\sqrt{x}+O(d^2).
\end{multline}

\smallskip

Enfin, en utilisant la relation asymptotique
$$
\fhi(t)^2=\frac 1{t^4} +O(1/t^5) \quad (t\soe 1)
$$
et la proposition \ref{t68}, on obtient 
\begin{multline}\label{t55}
-x^2\Big (\sum_{K_{d} < k\ioe K_{d+1}}-\sum_{K_{d-1}<k\ioe K_{d}}\Big)\fhi(x/k)^2 =\\
-\frac{1}5 \big((d+1)^2\sqrt{d+1}+(d-1)^2\sqrt{d-1}-2d^2\sqrt{d}\,\big)\sqrt{x}+O(d^2).
\end{multline}

\subsubsection{Estimation de $Q_d(x)$}

En ins\'erant \eqref{t71}, \eqref{t72}, \eqref{t73}, \eqref{t52}, \eqref{t53} et \eqref{t55} dans \eqref{t56},  on obtient la proposition suivante.
\begin{prop}\label{t57}
Pour $d \in \Nat^*$ et $x>0$, on a
$$
Q_d(x)=\vartheta(d)\sqrt{x}+O(d^2).
$$
o\`u
$$
\vartheta(d)=\frac{32d^2}5\sqrt{d}-\frac{48d^2+16d-2}{15}\sqrt{d+1}-\frac{48d^2-16d-2}{15}\sqrt{d-1}.
$$
\end{prop}
\dem

Il s'agit de voir que l'on peut omettre le terme d'erreur $O(d^{7/2}x^{-1/2})$ et la contrainte $d+1\ioe x/2$.

Si $d\ioe x^{1/3}$, ce terme d'erreur est absorb\'e par le $O(d^2)$. D'autre part, en utilisant l'approximation 
$$
\sqrt{d\pm 1}=\sqrt{d}\big(1\pm 1/2d-1/8d^2\pm 1/16d^{3}+O(d^{-4})\big)\, ,
$$
on voit que $\vartheta(d)=O(d^{-3/2})$. La majoration uniforme $Q_d(x)\ll \sqrt{x/d}$ montre alors
que
$$
Q_d(x)- \vartheta(d)\sqrt{x} \ll d \quad (d>x^{1/3}),
$$
ce qui permet encore d'omettre le terme d'erreur $O(d^{7/2}x^{-1/2})$. Enfin, l'\'enonc\'e de la proposition est trivial si $d>x/2-1$.\fin

\subsection{Sommation de la s\'erie des $\vartheta(d)$}

Puisque $\vartheta(d)=O(d^{-3/2})$, la s\'erie $\sum_{d\soe 1}\vartheta(d)$ converge ; nous allons calculer sa somme.

\smallskip

En \'ecrivant
\begin{multline*}
\vartheta(d)=\frac{16}5\big(d^{5/2}-(d+1)^{5/2}\big)- \frac{16}5\big((d-1)^{5/2}-d^{5/2}\big)\\
+\frac{16}3\big((d+1)^{3/2}-(d-1)^{3/2}\big)-2(d+1)^{1/2}-2(d-1)^{1/2}.
\end{multline*}
et et en supposant $D$ entier positif, on obtient comme dans \cite{bw2015}, \S 2.5,
\begin{align}
\sum_{d\ioe D}\vartheta(d) &=-\frac 2{15}-\frac{16}{5}\big((D+1)^{5/2}-D^{5/2}\big)\\
&\qquad +\frac{16}{3}\big((D+1)^{3/2}+D^{3/2}\big)-2\big((D+1)^{1/2}-D^{1/2}\big)-4\sum_{d=1}^D\sqrt{d}\notag\\
&=-\frac 2{15}+\frac{\zeta(3/2)}{\pi}+O(D^{-1/2}).\label{t36}
\end{align}

\subsection{Conclusion}

Pour $x>0$ et $D\soe 2$, on a
\begin{align*}
Q(x) &= Q_0(x)+\sum_{1\ioe d\ioe D}Q_d(x)+\sum_{d>D}Q_d(x)\\
&=\frac{2}{15}\sqrt{x}+O(1)+\sum_{1\ioe d\ioe D}\big(\vartheta(d)\sqrt{x}+O(d^2)\big) +O(\sqrt{x/D})\expli{d'apr\`es \eqref{t49}, la proposition \ref{t57}, et \eqref{t48}}\\
&=\big (2/15+\sum_{d=1}^{\infty} \vartheta(d)\big)\sqrt{x}+O(D^3)+O\big(\sqrt{x/D}\,\big)\expli{car $\vartheta(d)=O(d^{-3/2})$}\\
&=\frac{\zeta(3/2)}{\pi}\sqrt{x}+O(x^{3/7}),
\end{align*}
d'apr\`es \eqref{t36}, et en choisissant $D=x^{1/7}$.

\providecommand{\bysame}{\leavevmode ---\ }
\providecommand{\og}{``}
\providecommand{\fg}{''}
\providecommand{\smfandname}{et}
\providecommand{\smfedsname}{\'eds.}
\providecommand{\smfedname}{\'ed.}
\providecommand{\smfmastersthesisname}{M\'emoire}
\providecommand{\smfphdthesisname}{Th\`ese}

\footnotesize

\noindent BALAZARD, Michel\\
Aix Marseille Universit\'e, CNRS, Centrale Marseille, I2M UMR 7373\\
13453, Marseille\\
FRANCE\\
Adresse \'electronique : \texttt{balazard@math.cnrs.fr}

\end{document}